\documentclass[12pt]{article}
\usepackage{a4,amsmath,amssymb,amscd,cite,graphicx,multirow,slashbox,cancel}
\def\Bbb{\mathbb}
\def\BZ{\Bbb Z}

\def\Tr{\mathrm{Tr}}

\newcommand{\etabox}[2]{\underset{\ ~#2}{\mbox{\scriptsize $#1$}\ \framebox[15pt]{\phantom{a}}}}
\newcommand{\Mtwo}{M_{12}\!:\!2}
\catcode`@=11 \@addtoreset{equation}{section} \catcode`@=12

\begin{document}
\bibliographystyle{utphys}
\begin{titlepage}
\renewcommand{\thefootnote}{\fnsymbol{footnote}}
\noindent
{\tt IITM/PH/TH/2010/8}\hfill
\texttt{AEI-2010-171}\\[4pt]
\mbox{} \hfill{\texttt{[v1.0]\,Dec.\,2010}}

\begin{center}
\large{\sf  Brewing moonshine for Mathieu}
\end{center} 
\bigskip 
\begin{center}
{\sf Suresh Govindarajan}\footnote{\texttt{suresh@physics.iitm.ac.in}} \\[3pt]
\textit{\footnote{Permanent Address}Department of Physics, Indian Institute of Technology Madras,\\ Chennai 600036, India \\[4pt]
\textrm{and} \\[4pt]
Max-Planck-Institut f\"ur Gravitationsphysik, Albert-Einstein-Institut, \\
D-14424 Golm, Germany}
\end{center}
\bigskip
\bigskip
\centerline{To John McKay, for inspiration}
\bigskip
\bigskip
\begin{abstract}
We propose a moonshine for the sporadic Mathieu group $M_{12}$ that relates its conjugacy classes to various modular forms and Borcherds Kac-Moody superalgebras. 
\end{abstract}
\end{titlepage}
\setcounter{footnote}{0}
\section{Introduction}

\textit{\small Moonshine is not a well defined term, but everyone in the area recognizes it when
they see it. Roughly speaking, it means weird connections between modular forms and sporadic simple groups. It can also be extended to include related areas such as infinite dimensional Lie algebras or complex hyperbolic reflection groups. Also, it should only be applied to things that are weird and special: if there are an infinite number of examples of something, then it is not moonshine.} \hfill -- \textsf{\small R. E. Borcherds}\cite{Borcherds:2001} \\

In this paper, we propose a moonshine for the sporadic Mathieu group $M_{12}$ in the spirit of the above statement. This is best summarized in the following figure:
\begin{figure}[hbt]
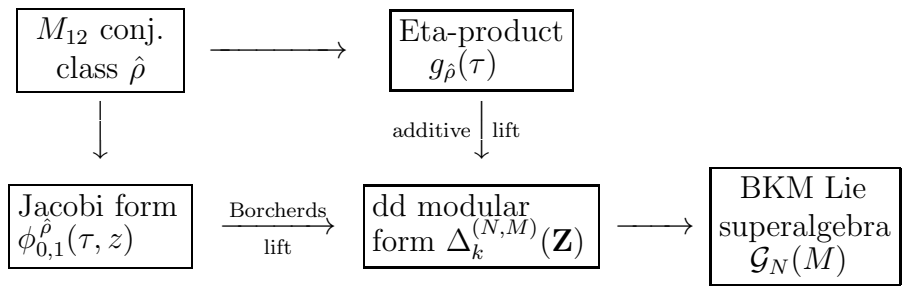

$$
\begin{CD}
\fbox{\parbox{2cm}{~$M_{12}$ conj. \\ \mbox{}~~~class $\hat{\rho}$}} @>\phantom{\textrm{Borcherdsab}}>> \fbox{\parbox{2.2cm}{Eta-product \\ \hspace*{0.4cm}$g_{\hat{\rho}}(\tau)$}} @. @.\\
@VVV  @V\textrm{additive}V\textrm{lift}V @. @. \\
\fbox{\parbox{2.2cm}{Jacobi form \\ $\phi^{\,\hat{\rho}}_{0,1}(\tau,z)$}} 
@>\textrm{Borcherds}>\textrm{lift}>  \fbox{\parbox{2.9cm}{dd modular\\ form $\Delta^{(N,M)}_k(\mathbf{Z})$}} @>>> \fbox{\parbox[c]{2.3cm}{~~BKM Lie \\
superalgebra $\mbox{~~~}\mathcal{G}_N(M)$}}
\end{CD}
$$
\caption{The proposed moonshine correspondence for $M_{12}$}\label{M12moonshine}
\end{figure}

Recent results  have provided evidence for the existence of a moonshine for the Mathieu group $M_{24}$\cite{Govindarajan:2009qt,Eguchi:2010ej,Cheng:2010pq,Gaberdiel:2010ch,Gaberdiel:2010ca,Eguchi:2010fg}. This moonshine in its most general form relates conjugacy classes of $M_{24}$ to genus-two Siegel modular forms that arise in the enumeration of dyonic degeneracies in a family of $\mathcal{N}=4$ string theories\cite{Dijkgraaf:1996it,Jatkar:2005bh,David:2006ud,David:2006ji} (see \cite{Sen:2007qy} for a review). Somewhat mysteriously, in some cases the \textit{square-root} of these Siegel modular forms appear as the Weyl-Kac-Borcherds (WKB) denominator formulae for Borcherds-Kac-Moody (BKM) Lie superalgebras\cite{Govindarajan:2008vi,Cheng:2008kt}.  This work is an attempt at understanding the `square-root' in terms of an outer automorphism of $M_{12}$ that leads to $M_{24}$. 

The canonical example that illustrates these ideas is the conjugacy class $1A (1^{24})$ of $M_{24}$. The $M_{24}$ moonshine maps this to the weight ten Igusa cusp form, $\Phi_{10}(\mathbf{Z})$. It's square-root is a weight five Siegel modular form (with character) $\Delta_5(\mathbf{Z})$. Gritsenko and Nikulin, in their studies of rank three Lorentzian Kac-Moody algebras\cite{Nikulin:1995,Gritsenko:2002},  have shown that $\Delta_5(\mathbf{Z})$ arises as the WKB denominator formula of a BKM Lie superalgebra. Clery and Gritsenko have constructed a family of modular forms that they call \textit{dd modular forms} that generalize $\Delta_5(\mathbf{Z})$\cite{Gritsenko:2008}. In an earlier paper\cite{Govindarajan:2010fu}, we have shown that all \textit{dd modular forms} appear as the square roots of Siegel modular forms that enumerate dyon degeneracies and that they arise as the WKB denominator formulae for rank three BKM Lie superalgebras. The \textit{main result} 
of this paper is to show that \textit{all} dd modular forms and their associated BKM Lie superalgebras are associated with a generalized moonshine for the sporadic group $M_{12}$.

The organization of the paper is as follows. In section two, we summarise the group theoretic aspects that are relevant for our considerations. We work out the relationship between $M_{12}$ and $M_{24}$ and show how one can track conjugacy classes of $M_{24}$ to those of $M_{12}$. In section three, we work out the first class of examples. These associate conjugacy classes of $M_{12}$ with balanced cycle shapes  to multiplicative eta-products, weight zero Jacobi forms of index one and dd modular forms. In section four, we consider a generalized moonshine for $M_{12}$ in the sense of Norton. This leads to the other dd modular forms as well as weight zero Jacobi forms of index $>1$. We conclude in section 5 with some remarks. Appendix A has some of the relevant group theoretic details while appendix B has some of the background material on modular forms.
\smallskip

\noindent \textbf{Notation:} We use a hat to distinguish objects associated with $M_{12}$ from those associated with $M_{24}$. Thus a conjugacy class of $M_{12}$ will be indicated by $\hat{\rho}$ while that of $M_{24}$ will be $\rho$. Characters of $M_{12}$ will be thus written as $\widehat{\chi}_i$ while characters of $\Mtwo$ are written with a tilde: $\widetilde{\chi}_i$.

\section{Group Theory}
\subsection{A quirk}

As described by Mark Ronan\cite{Ronanweb}, the Mathieu groups $M_{12}$ and $M_{24}$ arise due to the existence of certain quirks. Among all the permutation groups, only  $S_6$ admits an outer automorphism of order two that leads to $M_{12}$. In particular, one can show that the group $S_6\rtimes \BZ_2$ constructed using this outer automorphism is a maximal subgroup of $M_{12}$. Similarly,  $M_{12}$ admits an outer automorphism of order two that leads to $M_{24}$ in a similar fashion. Let us denote this automorphism of $M_{12}$ by $\varphi$ and the image of an element $g\in M_{12}$ under this automorphism by $\varphi(g)$. The group
$\Mtwo \equiv M_{12}\rtimes \mathbb{Z}_2$ is given by the set $M_{12} \times\BZ_2$ with the composition rule:
\begin{equation}
(g_1, h_1) \cdot (g_2,h_2) = (g_1 \cdot h_1(g_2), h_1\cdot h_2)\ ,
\end{equation}
where $g_1,g_2\in M_{12}$ and $h_1,h_2\in \BZ_2$ and $h(g)=g$ when $h=e$ and $h(g)=\varphi(g)$ when $h=\varphi$.

Now consider the realization of $M_{12}$ as a subgroup of the permutation group $S_{12}$ and let us use the same symbol $g$ to now indicate the $12\times 12$ permutation matrix in this realization. The $24$-dimensional representation of the group
$\Mtwo$ is then given by
\begin{equation}
(g,e) =  \begin{pmatrix} g & 0 \\ 0 & \varphi(g)\end{pmatrix}\quad, \quad
(g,\varphi) = \begin{pmatrix}  0 & g \\  \varphi(g) & 0\end{pmatrix}\quad \forall g\in M_{12}\ . \label{M12Mtwo}
\end{equation}
The group $\Mtwo$ is a maximal subgroup of $M_{24}$. In the sequel, the  conjugacy classes associated with elements of type $(g,e)\in \Mtwo$ will play an important role in our considerations.

\subsection{Conjugacy classes of $M_{12}$ and $\Mtwo$}

$M_{12}$ has fifteen conjugacy classes and the outer automorphism $\varphi$ acts on the conjugacy classes of $M_{12}$. It interchanges the conjugacy classes (cycle shapes),
$$
4A\ (2^24^2) \leftrightarrow 4B\ (1^44^2)\quad ; \quad  8A\ (4~8)\leftrightarrow 8B\ (1^22~8)\quad  \textrm{and}\quad  11A \leftrightarrow 11B\ ,
$$ 
leaving all other conjugacy classes invariant. This observation enables us to track how conjugacy classes of $M_{12}$ combine into conjugacy classes of $\Mtwo$ using the 24-dimensional representation that we just constructed. For instance, the cycle shape $1^82^8$($2B$) of $\Mtwo$ decomposes into two identical copies of the $M_{12}$ conjugacy class $1^{4}2^4$($2B$), the cycle shape $1^42^24^4$($4A$) of $\Mtwo$ decomposes as $1^44^2$($4B$) and $2^24^2$($4A$). It is easy to see that both these elements arise in $\Mtwo$ in the form of $(g,e)$. Of course, $\Mtwo$ has conjugacy classes that do \textit{not} reduce to conjugacy classes of $M_{12}$ in this fashion. One such class is the one corresponding to the cycle shape $4^6(4C)$.

Since both $\Mtwo$ and $M_{24}$ have 24 dimensional (permutation) representations, it is rather easy to  track conjugacy classes directly in terms of cycle shapes. We find that among the cycle shapes that appear in the half-BPS counting, only the cycle shape $1^37^3(7A)$ does \textit{not} appear as a conjugacy class of $\Mtwo$. This implies that all symplectic automorphisms of K3 other than the $\BZ_7$ one can also be considered as elements of $\Mtwo$.
In Table \ref{cycleshapes}, we track how  some conjugacy classes of $M_{24}$  
realized by group elements of the form $(g,e)$ given in Eq. \eqref{M12Mtwo} decompose into $M_{12}$ conjugacy classes.

\subsubsection*{Balanced cycle shapes}

 A cycle shape, $\rho=1^{a_1}2^{a_2}\cdots N^{a_N}$, is said to be balanced if there exists a positive integer $M$ such that $\big(\tfrac{M}1\big)^{a_1} \big(\tfrac{M}2\big)^{a_2}\cdots \big(\tfrac{M}n\big)^{a_n}$ is the same as $\rho$. It is known that all conjugacy classes of $M_{24}$ arise from balanced cycle shapes\cite{Conway:1979}. However, that is not necessarily true for all $M_{12}$ conjugacy classes. We observe that the cycle shapes associated with the $M_{12}$ conjugacy classes $4B$ and $8B$ are \textit{not} balanced.

\begin{table}[hbt]
\centering
\newcommand\T{\rule{0pt}{2.6ex}}
\begin{tabular}{ccccccccc} \hline
$\rho$ & $1^{24}$ &  $1^82^8$ & $1^63^6$ & $1^42^24^4$ & $1^45^4$ & $1^22^23^26^2$ & $1^37^3$ & $1^2 2^14^18^2$ \T \\[3pt] \hline
$M_{24}$  class & 1A & 2A & 3A & 4B & 5A & 6A & 7A & 8A \\[3pt]
$\Mtwo$  class & 1A & 2B & 3A & 4A & 5A & 6B & -- & 8A \\[3pt]
$M_{12}$ classes & 1A/1A & 2B/2B & 3A/3A & 4A/4B & 5A/5A & 6B/6B & -- & 8A/8B \\[3pt]
 \hline
\end{tabular}
\caption{From $M_{24}$ cycle shapes to $M_{12}$ cycle shapes\cite{Atlasv3}.}\label{cycleshapes}
\end{table}

\subsection{Irreps of $M_{12}$ and $\Mtwo$}

The decomposition of irreps of $\Mtwo$ into those of $M_{12}$ can also be worked out. For instance, $22=11\oplus 11'$ where $11'$ is the image of $11$ under the action of the outer automorphism $\varphi$. Similarly $16'$ is the image of $16$ and $55'$ is the image of $55$.  All other irreps of $M_{12}$ are invariant under the action of $\varphi$. Such irreps arise from two inequivalent irreps of $\Mtwo$ of the same dimension.

\section{The correspondence}

\subsection{Additive $M_{12}$ moonshine}

The additive version of the $M_{24}$ moonshine correspondence relates cycle shapes  of $M_{24}$ to multiplicative eta-products through the map\cite{Mason:1985,Dummit:1985,Booherthesis}:
\begin{equation}\label{cycleshapemap}
\rho=1^{a_1}2^{a_2}\cdots N^{a_N} \longmapsto g_\rho(\tau)\equiv \prod_{j=1}^N\eta(j\tau)^{a_j}\ ,
\end{equation}
where $\eta(\tau)$ is the Dedekind eta function.
These eta-products appear as the generating function of the degeneracy of twisted electrically charged $\tfrac12$-BPS states in type IIA string theory compactified on $K3\times T^2$\cite{Govindarajan:2009qt}. Among the conjugacy classes of $M_{24}$ that appear in this fashion, all classes other than  the $7A=1^37^3$ conjugacy class reduce to conjugacy classes of $\Mtwo$.

We \textit{propose} that  the additive version of the $M_{12}$ moonshine correspondence relates $M_{12}$ conjugacy classes to eta-products through the map Eq. \eqref{cycleshapemap}. Thus, the eta-product for the $M_{12}$ conjugacy class $2B$ is $\eta(\tau)^4\eta(2\tau)^4$. In all cases except for $4A$ and $8A$,\footnote{This happens as the outer-automorphism of $M_{12}$ acts non-trivially on these two conjugacy classes.} one observes that this corresponds to taking the square-root of a $M_{24}$ eta-product.  This is similar to McKay's observation (proved in \cite{Kac:1980}) that the cube-root of the modular invariant $j$ provides a moonshine for the group $E_8$. We thus have the relationship
\begin{equation}\label{spliteta}
\boxed{
g_\rho(\tau) = g_{\hat{\rho}}(\tau) \times g_{\varphi(\hat{\rho})}(\tau)
}\ ,
\end{equation}
where $\rho$ is a conjugacy class of $M_{24}$ (as well as $\Mtwo$) and $\hat{\rho}$ is a conjugacy class of $M_{12}$.
In the last row of Table \ref{cycleshapes},  we write $\hat{\rho}/\varphi(\hat{\rho})$ to indicate the two conjugacy classes of $M_{12}$ into which a given conjugacy class of $M_{24}$ decomposes.  

In appendix \ref{multeta}, we show that all $M_{12}$ conjugacy classes with balanced cycle shapes give rise to multiplicative eta-products. Unlike $M_{24}$ whose conjugacy classes all have balanced cycle shapes, two $M_{12}$ conjugacy classes $4B$ and $8B$ have unbalanced cycle shapes and hence we shall exclude them from most of our considerations. We shall also distinguish between $M_{12}$ conjugacy classes that reduce to conjugacy classes of $M_{11}$ and those that don't. The one's that do (reduce to $M_{11}$ conjugacy classes)  have cycle shapes with at least one one-cycle. 

The $M_{12}$-module that leads to the eta-product is easy to construct. It is given by the Fock space, $\mathcal{F}^\natural$,  of the oscillator (non-zero) modes of $12$ chiral scalars. As can be seen from a computation more or less identical to the one discussed in the appendix A of \cite{Govindarajan:2010fu}. The following trace over the Fock space can be expressed in terms of $M_{12}$ characters. One has 
\begin{align}
\frac1{g_{\hat{\rho}}(\tau)} &= \Tr_{\mathcal{F}^\natural}\Big( g\ q^{L_0-1/2}\Big)\ , \nonumber \\
&= q^{-1/2}\Big(1+ \big[1+\widehat{\chi}_2(\hat{\rho})\big]\  q + \big[3+3\widehat{\chi}_2(\hat{\rho})+\widehat{\chi}_7(\hat{\rho})\big]\  q^2+\cdots\Big)\ ,
\end{align}
where $g$ is  an element of $M_{12}$ (in the conjugacy class $\hat{\rho}$ with balanced cycle shape) acting as a subgroup of the permutation group on the 12 chiral scalars and the characters are numbered as in the $M_{12}$ character table(see Eq. \eqref{M12chartable}). 
One can also see the appearance of $M_{12}$ characters in the Fourier coefficients expansion of $q^{-1/2}\ g_{\hat{\rho}}(\tau)$. One also has
\begin{equation}
q^{-1/2}\ g_{\hat{\rho}}(\tau) = 1 -\big[1+\widehat{\chi}_2(\hat{\rho})\big]\  q + \big[-1 + \widehat{\chi}_9(\hat{\rho})\big]\ q^2 + \cdots \ .
\end{equation}

\subsection{Multiplicative $M_{12}$ moonshine}

As we just did for the additive $M_{12}$ moonshine, we shall first consider the multiplicative version of moonshine that appears for $M_{24}$\cite{Eguchi:2010ej,Cheng:2010pq,Gaberdiel:2010ch,Gaberdiel:2010ca,Eguchi:2010fg} and show that it implies a multiplicative moonshine for $M_{12}$ as well. The multipicative moonshine correspondence for $M_{24}$ maps its conjugacy classes to weight zero Jacobi forms that arise as twisted elliptic genera of $K3$\cite{David:2006ji,David:2006ud,David:2006yn}. The appearance of representations of $M_{24}$ is seen by decomposing the elliptic genus in terms of characters of the $\mathcal{N}=4$ superconformal algebra (SCA) at level $k=1$\cite{Eguchi:2008gc,Eguchi:2009cq,Eguchi:2009ux}. In the expansion below, $\alpha$ is a constant while $\Sigma(\tau)$ is a function of $\tau$. 
\begin{equation}
\psi^\rho_{0,1}(\tau,z) = \alpha^\rho~ \mathcal{C}(\tau,z) + q^{-\tfrac18}~\Sigma^\rho(\tau)~ \mathcal{B}(\tau,z)\ ,
\end{equation}
where $\mathcal{C}(\tau,z)$ is the massless character and $\mathcal{B}(\tau,z)$ is the massive character of the $\mathcal{N}=4$ SCA  at level one and $\rho$ denotes a $M_{24}$ conjugacy class. The characters are such that 
$$
\mathcal{C}(\tau,z=0)=1\quad \textrm{ and }\quad \mathcal{B}(\tau,z=0)=0\ .
$$
Hence, it is easy to show that $\alpha^\rho = \chi^\rho(K3)$, the twisted Euler characteristic of K3.  The $M_{24}$ moonshine correspondence implies that 
$$
\alpha^\rho=\chi_1(\rho)+\chi_{23}(\rho)\ ,
$$ 
where $\chi_{N}$ represents the character for the $M_{24}$-irrep of dimension $N$\cite{Eguchi:2010ej,Cheng:2010pq}.
The function $\Sigma^\rho(\tau)$ has the following Fourier expansion
\begin{equation}
\Sigma^\rho(\tau)=\Big(-2 + \sum_{n=1}^\infty A^\rho(n)~ q^{n}\Big)\ ,
\end{equation}
where  $A^\rho(n)$ are integers that can also be expressed in terms of characters of $M_{24}$\cite{Eguchi:2010ej,Cheng:2010pq,Gaberdiel:2010ch,Gaberdiel:2010ca,Eguchi:2010fg}. One has
\begin{multline}
\Sigma^\rho(\tau)= -2 + [\chi_{45}(\rho)+ \chi_{\overline{45}}(\rho) ]\ q + [\chi_{231}(\rho) +\chi_{\overline{231}}(\rho) ]\ q^2 \\ + [\chi_{770}(\rho)+ \chi_{\overline{770}}(\rho) ]\ q^3   + [\chi_{2277}(\rho)+ \chi_{\overline{2277}}(\rho) ]\ q^4 + \cdots   
\end{multline}

As we did for the case of the additive moonshine, we look to rewrite the conjugacy classes of  $M_{24}$ that appear in the multiplicative moonshine in terms of $M_{12}$ conjugacy classes. In other words, we seek a relationship of the form
\begin{align}
\alpha^\rho &= \widehat{\alpha}^{\,\hat{\rho}} + \widehat{\alpha}^{\, \varphi(\hat{\rho})}\ , \\
\Sigma^\rho(\tau) & = \widehat{\Sigma}^{\,\hat{\rho}}(\tau) +\widehat{\Sigma}^{\,\varphi(\hat{\rho})}(\tau)\ ,
\end{align}
where we decompose the $M_{24}$ conjugacy class $\rho$ in terms of two conjugacy classes of $M_{12}$ that we denote by $\hat{\rho}$ and its image under the outer automorphism $\varphi(\hat{\rho})$. This works for all $M_{24}$ conjugacy classes given in Table except for $7A$ which does not reduce to a conjugacy class of $\Mtwo$.

Given a conjugacy class $\hat{\rho}$ of $M_{12}$, we consider the Jacobi form of weight zero and index one given by
\begin{equation}\label{characterexpand}
\phi^{\hat{\rho}}_{0,1}(\tau,z) = \widehat{\alpha}^{\,\hat{\rho}}~ \mathcal{C}(\tau,z) + q^{-\tfrac18}~\widehat{\Sigma}^{\,\hat{\rho}}(\tau)~ \mathcal{B}(\tau,z)\ ,
\end{equation}
where  $\widehat{\alpha}^{\,\hat{\rho}}= 1+\widehat{\chi}_2(\hat{\rho})$ and
\begin{multline}
\widehat{\Sigma}^{\hat{\rho}}(\tau)= -1 + \widehat{\chi}_6(\hat{\rho})\ q + [\widehat{\chi}_{8}(\hat{\rho})+\widehat{\chi}_{15}(\hat{\rho}) ]\ q^2  \\ + [\widehat{\chi}_{11}(\hat{\rho})+2\ \widehat{\chi}_{13}(\hat{\rho})+2\ \widehat{\chi}_{14}(\hat{\rho})+\widehat{\chi}_{15}(\hat{\rho}) ]\ q^3   + \cdots   \ .
\end{multline}
The above formulae are obtained by using the decomposition of $M_{24}$ irreps into $M_{12}$ irreps as given in appendix \ref{decompose}. Thus, one has
\begin{equation}\label{JacobiM12split}
\boxed{
\psi_{0,1}^\rho(\tau,z) = \widehat{\psi}_{0,1}^{\,\hat{\rho}}(\tau,z) +\widehat{\psi}_{0,1}^{\,\varphi(\hat{\rho})}(\tau,z)
}\quad .
\end{equation}
It is easy to see that for conjugacy classes that are invariant under the outer automorphism, one has $\psi^\rho_{0,1}(\tau,z)=2\widehat{\psi}^{\hat{\rho}}_{0,1}(\tau,z)$. However,  the $M_{24}$ conjugacy classes $\rho=4B,8A$ decompose into  distinct $M_{12}$ conjugacy classes, one of which is not balanced and we will not consider them.

\subsection{The genus-two Siegel modular form}

\subsubsection{The additive lift}

Let $g$ denote an element of $M_{24}$ (that reduces to an element of $M_{23}$) whose conjugacy class is $\rho$. It has been shown in ref. \cite{Govindarajan:2010fu} that the degeneracy of $g$-twisted $\tfrac14$-BPS multiplets in the $T^6$ compactified heterotic string is generated by a genus-two Siegel modular form. This Siegel modular form is obtained as the additive lift of a Jacobi form of weight $k$, index $1$ and level $N$ given by
\begin{equation}
\phi^\rho_{k,1}(\tau,z) \equiv   \frac{\theta_1(\tau,z)^2}{\eta(\tau)^6} \times 
g_\rho(\tau)\ ,
\end{equation}
where $k=\tfrac12(\sum_i a_i)-2$ and we have indicated the $M_{24}$ conjugacy class as a superscript. The Siegel modular form is then given by the additive lift of the Jacobi form, $\phi_{k,1}(\tau,z)$. One has\cite{Jatkar:2005bh,Govindarajan:2010fu}
\begin{equation}
\Phi^{\rho}_k(\mathbf{Z}) = \mathcal{A}\Big[\phi_{k,1}(\tau,z)\Big] = \sum_{m=1}^\infty \phi^\rho_{k,m}(\tau,z) ~ s^m\ ,
\end{equation}
where $\mathbf{Z} \equiv \left(\begin{smallmatrix} \tau & z \\ z & \sigma \end{smallmatrix}\right)\in \mathbb{H}_2$, $r=\exp(2\pi i z)$ and $\phi^\rho_{k,m}(\tau,z)$ is a Jacobi form of weight $k$ and index $m$ obtained by the action of the Hecke operator on the additive seed $\phi^\rho_{k,1}(\tau,z)$:
\begin{align}
\phi^\rho_{k,m}(\tau,z)  &\equiv T_-^{(N)}(m) \phi_{k,1}(\tau,z)  \nonumber \\
&= \frac1m  \sum_{\substack{ad=m\\ (a,N)=1 }} \sum_{b=0}^{d-1}\chi(a)\ a^{k} \
\phi^\rho_{k,1}(\tfrac{a\tau+b}{d}, a z)  \ . \label{DMVVformula}
\end{align}
This is \textit{not} the most general form of the additive lift. We have given the simplest case in order to emphasize the fact that the new Jacobi forms are determined completely in terms of the additive seed\footnote{The most general Hecke operator that appears in the $\tfrac14$-BPS counting is discussed by Clery and Gritsenko\cite{Gritsenko:2008,Govindarajan:2010fu}.}.

Given that the eta-product uniquely determines the Siegel modular form, $\Phi^\rho_k(\mathbf{Z})$ one anticipates that this Siegel modular form should also be obtained as a trace over some module graded by $M_{24}$. Evidence towards the veracity of this statement is provided in ref. \cite{Govindarajan:2010cc}. Given our observation in Eq. \eqref{spliteta},  that the $M_{24}$ eta-products split into two $M_{12}$ eta-products, we expect something similar to happen with the Siegel modular form, $\Phi^\rho_k(\mathbf{Z})$. In other words, we  expect in all the cases where the eta product split, the Siegel modular form also splits into the product of two other Siegel modular forms  as follows:
\begin{equation}\label{splitsiegel}
\boxed{
\Phi^\rho(\mathbf{Z}) = \Delta^{\hat{\rho}}(\mathbf{Z}) \times \Delta^{\varphi(\hat{\rho})}(\mathbf{Z})\ ,
}
\end{equation}
where $\rho$ is a conjugacy class of $M_{24}$ (as well as $\Mtwo$) and $\hat{\rho}$ is a conjugacy class of $M_{12}$. Let us focus on conjugacy classes of $M_{12}$ that are invariant under the outer automorphism. In particular, consider the conjugacy classes $1A$, $2B$, $3A$ and $5A$.  In all these cases, this implies that
\begin{equation}
\Phi^\rho_k(\mathbf{Z}) = \Big[\Delta_{k/2}^{\hat{\rho}}(\mathbf{Z})\Big]^2
\end{equation}
This indeed agrees with the observations made in ref. \cite{Govindarajan:2008vi} where it was shown that the Siegel modular forms $\Delta_{k/2}^{\hat{\rho}}(\mathbf{Z})$ are natural generalizations of $\Delta_{5}(\mathbf{Z})$ which is the square-root of the weight ten Igusa cusp form. Further, it was shown that they are given by the additive lift
\begin{equation}\label{M12additivelift}
\Delta^{\hat{\rho}}_{k/2}(\mathbf{Z}) = \mathcal{A}\Big[\tfrac{\theta_1(\tau,z)}{\eta(\tau)^3} \times 
g_{\hat{\rho}}(\tau)
\Big] \ ,
\end{equation}
where $g_{\hat{\rho}}(\tau)$ is the eta-product associated with the $M_{12}$ conjugacy class, $\hat{\rho}$. In other words, the additive lift given above provides a direct link between the Siegel modular forms $\Delta_{k/2}(\mathbf{Z})$ and $M_{12}$ conjugacy classes that parallels what happened between the $\Phi_k(\mathbf{Z})$ and $M_{24}$ conjugacy classes. In most examples, these Siegel modular forms are examples of the \textit{dd modular forms} of Clery and Gritsenko\cite{Gritsenko:2008}.

\subsubsection*{The splitting of conjugacy class $4B$} 

The $M_{24}$ conjugacy class $4B$ is an interesting one. From table \ref{cycleshapes}, we see that it splits into two \textit{distinct} $M_{12}$ conjugacy classes $4A$ and $4B$.\footnote{This splitting is different from the square-root of $\Phi^{4B}_3(\mathbf{Z})$ considered in \cite{Govindarajan:2009qt} The additive seed in that case arises from the cycle shape $1^224^2$ which does not occur as a $M_{12}$ conjugacy class. 
}. This suggests that we should expect a splitting of the form:
$$
\Phi^{4B}_3(\mathbf{Z}) = \Delta_2^{4A}(\mathbf{Z})\times  \Delta_1^{4B}(\mathbf{Z}) \ ,
$$
where the modular forms  $\Delta^{\hat{\rho}}(\mathbf{Z})$ are given by the additive lift given in Eq. \eqref{M12additivelift}. In particular, the Siegel modular form, $\Delta_1^{4A}(\mathbf{Z})$,  corresponding to $\hat{\rho}=4A$ (cycle shape $2^24^4$) has already appeared in \cite[see sec.~3]{Govindarajan:2010fu} where it was denoted by $Q_1(\mathbf{Z})$ and can be expressed in terms of genus-two theta constants (see appendix B). The Siegel modular form $\Phi^{4B}_3(\mathbf{Z})$ was constructed in \cite{Govindarajan:2009qt} and is  also expressible in terms of products of genus-two theta constants. The modular form for $\Delta_2^{4B}(\mathbf{Z})$ has not been constructed so far. However, it can be expressed as the quotient of $\Phi^{4B}_3(\mathbf{Z})$ by $\Delta_2^{4A}(\mathbf{Z})$. 
\begin{equation}
\Delta_2^{4B}(\mathbf{Z})= \frac{\Phi^{4B}_3(\mathbf{Z})}{\Delta_1^{4A}(\mathbf{Z})}   \ .
\end{equation}
Does it imply that it is a meromorphic form? Interestingly, that doesn't happen as all the genus-two theta constants that appear in $\Delta_1^{4A}(\mathbf{Z})$ also occur in $\Phi^{4B}_3(\mathbf{Z})$ leading to a neat cancellation! Thus, $\Delta_2^{4B}(\mathbf{Z})$ is \textit{not} a meromorphic modular form. Further, the $M_{12}$ conjugacy class $4B$ does not have a balanced cycle shape and hence is not considered.

\subsubsection*{The conjugacy classes $6A$ and $8A$}

The $M_{24}$ conjugacy class $6A$ splits into two copies of the $M_{12}$ conjugacy class $6B$. Thus we expect it behave in a fashion similar to $1A$. The $M_{24}$ conjugacy class $8A$ behaves like the $4B$ conjugacy class and splits into two distinct conjugacy classes $8A/8B$ as can be seen from Table \ref{cycleshapes}. The corresponding modular forms have been constructed by the additive lift in \cite{Govindarajan:2009qt} but have not been extensively studied.  Further, the $M_{12}$ conjugacy class $8B$ does not have a balanced cycle shape.

\subsubsection{The multiplicative lift}

We have seen that the twisted elliptic genera of $K3$, $\psi^\rho_{0,1}(\tau,z)$,  which provide a multiplicative moonshine for $M_{24}$ also can be written as the sum of two terms each arising from a conjugacy class of $M_{12}$ as given in Eq. \eqref{JacobiM12split}. It turns out that these twisted elliptic genera lead to product formulae for the Siegel modular forms $\Phi_k^\rho(\mathbf{Z})$. Similarly, the Jacobi forms $\widehat{\psi}_{0,1}^{\hat{\rho}}(\tau,z)$ provide a  product formula for $\Delta^{\hat{\rho}}(\mathbf{Z})$.
In particular, for the conjugacy classes, $\rho=1A,2A,3A,5A$, of $M_{24}$, the relationship  $\psi^\rho_{0,1}(\tau,z)=2\ \widehat{\psi}_{0,1}^{\hat{\rho}}(\tau,z)$ is consistent with the $\Phi_k^\rho(\mathbf{Z})$ being the square of $\Delta^{\hat{\rho}}(\mathbf{Z})$. This should also hold in principle for the conjugacy class $6A$ but we have  not checked this as the multiplicative lift has not been constructed (see \cite{Govindarajan:2009qt,Sen:2009md} for some details in this regard).

However, the conjugacy class $4B$ appears to lead to a different Siegel modular form. In particular, the product formula for $\Delta_1^{4B}(\mathbf{Z})$ arises from 
a Jacobi form of weight zero and index \textit{two}\cite{Govindarajan:2010fu}. Thus, it appears to us that $\phi_{0,1}^{4B}(\tau,z)$ does not generate a product formula for the Siegel modular form generated by the additive lift $\Delta_1^{4B}(\mathbf{Z})$. As we will discuss later, $4B$ is a conjugacy class of $M_{12}$ that does not descend to a conjugacy class of $M_{11}$. Such conjugacy classes arise from a generalized $M_{12}$ moonshine in the sense of Norton that we will discuss in a later section.

\subsection{Borcherds-Kac-Moody algebras}

A very nice result is that the $M_{12}$ modular forms $\Delta^{\hat{\rho}}(\mathbf{Z})$ discussed in the previous section (for $\hat{\rho}=1A,2B,3A,5A$ ) arise as the Weyl-Kac-Borcherds denominator formula for a family of rank-three Lorentzian Kac-Moody superalgebras\cite{MR1438983,Govindarajan:2008vi,Govindarajan:2009qt}. All the BKM Lie superalgebras have identical simple real roots with  Cartan matrix
$$
\begin{pmatrix} 2 & -2 & -2 \\ -2 & 2 & -2 \\ -2 & -2 & 2 \end{pmatrix}\ ,
$$
but differ in their imaginary roots. In other words, we end up with the sequence summarized in figure 1 that takes us from conjugacy classes of $M_{12}$ to Borcherds-Kac-Moody Lie superalgebras. In particular we obtain

\begin{table}[hb]
\centering
\newcommand\T{\rule{0pt}{2.6ex}}
\begin{tabular}{ccccc}\hline
$M_{12}$ conj. class & $1A$ & $2B$ & $3A$  & $5A$\T \\[3pt] \hline
BKM Lie algebra & $\mathcal{G}_1(1)$ & $\mathcal{G}_1(2)$ &  $\mathcal{G}_1(3)$ &  $\mathcal{G}_1(5)$\T \\[3pt] \hline
\end{tabular}
\caption{Relating $M_{12}$ conjugacy classes to Lie superalgebras(notation of \cite{Govindarajan:2010fu}) }
\end{table}
\section{Generalized $M_{12}$ moonshine}

In our approach, the  conjugacy classes of $M_{12}$ that don't reduce to conjugacy classes of $M_{11}$ such as $4B$ naturally appear when one considers a generalized  moonshine in the sense of Norton\cite[see appendix by Norton]{Mason:1987}. Generalized moonshine is best described using notation that is standard in Conformal Field Theory(CFT). The character of a module, $\mathcal{H}$  in  CFT is  given by 
\begin{equation}
\etabox11\equiv \Tr_{\mathcal{H}} \Big(q^{L_0-\frac{c}{24}}\Big)\ ,
\end{equation}
where the box notation will be explained soon. Now let $g$ and $h$ denote \textit{commuting} symmetries of finite order of the CFT. Let $\mathcal{H}_h$ denote the $h$-twisted module in the orbifold of the original CFT by the group generated by $h$. Then, we define
\begin{equation}
\etabox1h\equiv \Tr_{\mathcal{H}_h} \Big(q^{L_0-\frac{c}{24}}\Big)\ ,
\end{equation}
and similarly, one might consider a more general situation with the insertion of $g$ in a trace over the $h$-twisted module.
\begin{equation}
\etabox{g}h\equiv \Tr_{\mathcal{H}_h} \Big(g\ q^{L_0-\frac{c}{24}}\Big)\ .
\end{equation}
Using these ideas from CFT, we use the same pictures to represent suitable traces over twisted modules though we don't always specify the details of the module.

The three different moonshines that we have discussed involving eta-products, Jacobi forms and Siegel modular forms all arise from taking $g$ to be an element of $M_{12}$  and taking the trace over a suitable module graded by $g$ and thus are of type $\etabox{g}1$. We shall discuss the generalized moonshine for each of these modular forms.

\subsection{Eta-products}

The module here is furnished by the oscillator Fock space of twelve chiral bosons -- this was denoted earlier by $\mathcal{F}^\natural$. Then, as we have already seen, one has
\begin{equation}
g_{\hat{\rho}}(\tau) = \Tr_{\mathcal{F}^\natural} \Big(g\ q^{L_0-\frac12}\Big)\quad \longleftrightarrow \quad \etabox{g}1\ ,
\end{equation} 
where $\hat{\rho}$ is the conjugacy class of $g$. Let $h$ denote an element of $M_{12}$ (of order $N$) that acts by permuting the $12$ chiral bosons. and $\mathcal{F}_h$ denote the $h$-twisted Fock space. Then, it is natural to consider the generalized moonshine of type $\etabox1h$. For all conjugacy classes with balanced cycle shapes, a calculation analogous to the one in \cite[see appendix A]{Govindarajan:2010fu} gives rise to the eta product with modified argument\footnote{The conjugacy classes $4B$ and $8B$ do not fit this and are not considered.}
\begin{equation}
g_{\hat{\rho}}\big(\tfrac{\tau}{N}\big) = \Tr_{\mathcal{F}_h} \Big(\ q^{L_0-\frac1{2N}}\Big)\quad \longleftrightarrow \quad \etabox{1}h\ ,
\end{equation} 
where $\hat{\rho}=[h]$. These turn out be square-roots of the eta-products that count $\tfrac12$-BPS states in the  $\BZ_N$ CHL orbifold -- the $\BZ_N$ being generated by an element of $\Mtwo$ of the form $(h,e)$.

Next, one can consider the more general case of two different commuting elements of $M_{12}$ and a generalized moonshine of type $\etabox{g}h$ . Again, one obtains an eta-product (let $g$ have order $M$ and $h$ has order $N$ as before)
\begin{equation}
g_{\hat{\rho}}\big(\tfrac{\tau}{N}\big) = \Tr_{\mathcal{F}_h} \Big(g\ q^{L_0-\frac1{2N}}\Big)\quad \longleftrightarrow \quad \etabox{g}h\ ,
\end{equation} 
Again, using the relationship that we have observed with $g$-twisted $\tfrac12$-BPS states in the CHL $\BZ_N$-orbifold, we obtain three conjugacy classes corresponding to the values: $(M,N)=(2,2),(2,4),(4,2),(3,3)$ -- the first one gives rise to the conjugacy class $2A$, the next two correspond to the conjugacy class $4A$ and the last one gives rise to the conjugacy class $4A$.

\begin{table}[hbt]
\centering
\newcommand\T{\rule{0pt}{2.6ex}}
\begin{tabular}{cccc} \hline
$\hat{\rho}$ & $2^{6}$ &  $2^24^2$ & $3^4$  \T \\[3pt] \hline
$M_{12}$ class & 2A & 4A & 3B   \\[3pt]
 \hline
\end{tabular}
\end{table}

\subsection{The Siegel modular forms}

Given a multiplicative eta-product associated with a generalized moonshine of arbitrary type, we can immediately construct a  Siegel modular form by the additive lift. This enables us to provide candidate Siegel modular forms for a generalized moonshine of similar type. 

The Siegel modular forms associated with moonshine of type $\etabox1h$ 
are obtained by writing  the modular forms that count $\tfrac14$-BPS states in the CHL orbifold generated by the element $(h,e)$ (of order $N$) as in Eq. \eqref{splitsiegel}. For the $M_{12}$ conjugacy classes $1A$, $2B$, $3A$, this implies that the Siegel modular form is the square-root of the modular form counting $\tfrac14$-BPS dyons. The modular forms that, in principle, lead to the conjugacy classes $5A$ and $6B$ have not been constructed and hence our proposal remains conjectural in these cases.

A similar approach leads to generalized moonshine of type $\etabox{g}h$ with $g$ (of order $M$) and $h$ (of order $N$) commuting elements of $M_{12}$. Now, the modular forms that one obtains are square-roots of the Siegel modular forms that count $g$-twisted $\tfrac14$-BPS states in the CHL $\BZ_N$ orbifold. This leads to Siegel modular forms for the conjugacy classes $2A$,  $3B$,  $4A$. The conjugacy class $4A$ has two possibilities corresponding to   $(N,M)=(2,4)$ and $(4,2)$. All these modular forms are again instances of the dd-modular forms of Clery and Gritsenko\cite{Gritsenko:2008}.

As with the other modular forms, these dd modular forms also arise as the Weyl-Kac-Borcherds denominator formulae for BKM Lie superalgebras considered in \cite{Govindarajan:2010fu}. We list them in Table \ref{bkmtable}. 
\begin{table}[hb]
\centering
\newcommand\T{\rule{0pt}{2.6ex}}
\begin{tabular}{c|cc|cccc}\hline
$M_{12}$ conj. class  & $2B$ & $3A$ & $2A$ & $3B$ & $4A$ &  $4A'$ \T \\[3pt] \hline
Modular Form &  $\Delta_3^{(2,1)}$ & $\Delta_2^{(3,1)}$ & $\Delta_2^{(2,2)}$ & $\Delta_1^{(3,3)}$ & $\Delta_1^{(2,4)}$ & $\Delta_1^{(4,2)}$\T \\
BKM Lie algebra  & $\mathcal{G}_2(1)$ &  $\mathcal{G}_3(1)$ &$\mathcal{G}_2(2)$ & $\mathcal{G}_3(3)$ & $\mathcal{G}_2(4)$ &$\mathcal{G}_4(2)$ \T \\[3pt] \hline
\end{tabular}
\caption{Siegel  modular forms and Lie superalgebras(notation of \cite{Govindarajan:2010fu}) for generalized moonshine } \label{bkmtable}
\end{table}

\subsection{Jacobi forms}

There are two routes to obtaining Jacobi forms associated with generalized moonshine.  The first method is look directly for generalization of the Jacobi forms of type $\etabox{g}1$ that we considered earlier. The second method is to consider the Jacobi forms that provide Borcherds/multiplicative lifts for the Siegel modular forms (of the previous subsection) associated with generalized moonshine. It turns out that these two methods do \textit{not} necessarily lead to the same Jacobi forms.

We begin with the first method and look
for  a  generalized moonshine for the Jacobi forms of type $\etabox1h$. The simplest way to obtain these Jacobi forms  is to consider their transformation under the $(\tau,z) \rightarrow (-1/\tau,z/\tau)$ of the Jacobi forms of type $\etabox{h}1$. In these examples,   the Jacobi forms that generate product formulae for the Siegel modular forms also lead to the same Jacobi form.

Let us first consider the conjugacy classes of $M_{12}$ that reduce to $M_{11}$ classes. In particular, consider the classes $1A$, $2B$, $3A$ and $5A$.  One expands the transformed Jacobi forms (of weight zero and index one) in terms of $\mathcal{N}=4$ characters as in Eq. \eqref{characterexpand}. However,  one can see that the function $\widehat{\Sigma}^{\hat{\rho}}(\tau)$ must have the following Fourier expansion
\begin{equation}
\widehat{\Sigma}^{\hat{\rho}}=\Big(\beta + \sum_{n=1}^\infty A^{\hat{\rho}}(n)~ q^{n/N}\Big)\ ,
\end{equation}
where the fractional power of $q$ reflects the width of the cusp and $A^{\hat{\rho}}(n)$ are \textit{conjectured} to be integers.  The coefficients $A(n)$ are conjectured to be integers to all orders. We find
\begin{align*}
\widehat{\Sigma}^{2B}(\tau) &= (8 q^{1/2} + 24 q + 56 q^{3/2} + 112 q^2 + \cdots ) \\
\widehat{\Sigma}^{3A}(\tau)&= (3 q^{1/3}+ 9 q^{2/3} + 15 q + 30 q^{4/3}+ 45 q^{5/3} +\cdots) \\
\widehat{\Sigma}^{5A}(\tau) &= (q^{1/5} + 3 q^{2/5} + 4 q^{3/5}+ 7q^{4/5}+ 9 q +\cdots)
\end{align*}

The next generalization is to consider generalized moonshine of type $\etabox{g}h$, where $g$ and $h$ are two commuting elements of $M_{12}$. From our multiplicative eta-products, we know that these lead to conjugacy classes $2A$, $3A$ and $4A$ -- these conjugacy classes do not reduce to conjugacy classes of $M_{11}$. The product formula for the corresponding Siegel modular forms have been discussed in \cite{Govindarajan:2010fu} and the associated Jacobi forms are of weight zero and index \textit{two}. Thus, they are somewhat different from the  other examples that we have considered. We can decompose these weight two Jacobi forms in terms of $\mathcal{N}=4$ superconformal characters. This will be discussed elsewhere\cite{Govindarajan:2010cc}.

In principle, we could consider the index one Jacobi forms associated with conjugacy classes of $M_{24}$ that do \textit{not} reduce to conjugacy classes of $M_{23}$. These were considered by Cheng as well as Gaberdiel et. al.\cite{Cheng:2010pq,Gaberdiel:2010ch,Gaberdiel:2010ca}. However, they do not seem to be related to Siegel modular forms to the best of our knowledge. Hence, we do not pursue this any further.

\section{Concluding Remarks}

In this paper, we have conjectured a relationship between $M_{12}$ conjugacy classes with balanced cycle shapes and  Siegel modular forms as well as BKM Lie superalgebras. The simplest examples appear for conjugacy classes of $M_{12}$ that reduce to conjugacy classes of $M_{11}$. We have proposed that other conjugacy classes correspond to a generalized moonshine. Evidence for these conjectured relationships have been provided for a large class of conjugacy classes -- notably for those classes with cycles less than $6$. In particular, we believe that the correspondence should hold for the $M_{12}$ conjugacy classes $6B$. We anticipate that there exists a Siegel modular form as well as BKM Lie superalgebra(s) associated with this conjugacy class\footnote{We thank Fabien Clery for an extensive email discussion in this regard.}\cite{Govindarajan:2011dd}. However, we are unsure about their existence for the conjugacy classes $6A$, $8A$, $10A$ and $11A/B$ -- the additive lift, if it exists, implies that these are modular functions (i.e., of weight zero).

We also anticipate that there exists a $M_{12}$-module $V^\natural$ that is graded by three integers $(n,\ell,m)$
$$
V^\natural =\oplus_{(n,\ell,m)}V_{(n,\ell,m)}\ ,
$$ 
 such that
\begin{equation}
\frac1{\Delta_5(\mathbf{Z})}= \sum_{(n,\ell,m)}\Tr_{V_{(n,\ell,m)}} \big(q^n r^\ell s^m\big)
\end{equation}
Insertions of elements $g\in M_{12}$  in the trace should lead to $\Delta_5(\mathbf{Z})$ being replaced by dd modular forms associated with the conjugacy class of $g$. Similarly, by considering $h$-twisted versions of the module $V^\natural$, we should recover all the dd modular forms.

In a forthcoming paper\cite{Govindarajan:2010cc}, we show that the Siegel modular forms $\Phi_k(\mathbf{Z})$ and $\Delta_k(\mathbf{Z})$ imply an infinite number of moonshines for the Mathieu groups $M_{24}$ and $M_{12}$ respectively. We also show that these moonshines include the additive and multiplicative moonshines that were discussed in this paper and elsewhere in the context of $M_{24}$.

The original motivation for our study of the $M_{12}$ moonshine was to understand the square-root that appeared in relating Siegel modular forms to BKM Lie superalgebras. However, our study has raised more questions than have been answered. We conclude with a few of these questions. Can we understand the appearance of the Fock space of 12 chiral bosons? Is there any relation to the worldvolume theory of a M5-brane wrapping a Enriques surface? Can we derive the dd modular forms in terms of a theory of multiple M5-branes?
\bigskip

\noindent \textbf{Acknowledgements:} A significant part of this work was done during a visit last summer at the Albert Einstein Institute at Potsdam. We thank all the members of the Institute and in particular, Stefan Theisen for a very productive stay. We also thank Karthik Inbasekar and  Dileep Jatkar for comments on a preliminary version of this paper as well as Prof. Naresh Dadhich for encouragement.

\appendix

\section{Character Tables}

\subsection{Character Table for $M_{12}$}\label{M12ch}

The character table for $M_{12}$ (obtained from the GAP character table database\cite{GAP4})
\begin{equation}
\left(\begin{smallmatrix}
\textrm{Label}         &1A &  2A & 2B & 3A &  3B &  4A & 4B & 5A & 6A & 6B & 8A & 8B & 10A &  11A & 11B \\[4pt]
\widehat{\chi}_1      &~1  &~1  &~1  &~1  &~1  &~1  &~1  &~1  &~1  &~1  &~1  &~1   &~1   &~1   &~1 \\
\widehat{\chi}_2   &   11& -1&  ~3&  ~2 &-1 &-1 & ~3 & 1 &-1  &~0 &-1 & ~1 & -1   &~0   &~0 \\
\widehat{\chi}_3    &  11'& -1 & ~3 & ~2& -1&  ~3& -1& ~1& -1  &~0  &~1 &-1 & -1   &~0   &~0 \\
\widehat{\chi}_4    &  16&  ~4  &~0& -2  &~1  &~0  &~0 & ~1 & ~1  &~0  &~0  &~0 & -1  & ~\alpha & \alpha^* \\
\widehat{\chi}_5    &  16' & ~4  &~0 &-2&  ~1  &~0  &~0  &~1 & ~1  &~0  &~0  &~0  &-1 & \alpha^* &  ~\alpha\\
\widehat{\chi}_6    &  45 & ~5 &-3  &~0 & ~3&  ~1 & ~1  &~0 &-1  &~0 &-1 &-1   &~0 &  ~1&   ~1\\
\widehat{\chi}_7  &    54 & ~6  &~6  &~0  &~0 & ~2 & ~2& -1  &~0  &~0  &~0  &~0  & ~1  &-1  &-1\\
\widehat{\chi}_8  &    55_R &-5  &~7  &~1 & ~1 &-1& -1  &~0 & ~1 & ~1 &-1& -1   &~0   &~0   &~0 \\ 
\widehat{\chi}_9  &    55& -5 &-1  &~1 & ~1&  ~3& -1  &~0 & ~1 &-1 &-1 & ~1   &~0   &~0   &~0 \\ 
\widehat{\chi}_{10}&     55'& -5 &-1&  ~1 & ~1 &-1 & ~3  &~0 & ~1& -1&  ~1 &-1   &~0   &~0   &~0\\
\widehat{\chi}_{11} &   66 & ~6  &~2&  ~3  &~0 &-2& -2&  ~1  &~0& -1  &~0  &~0 & ~1   &~0   &~0\\
\widehat{\chi}_{12} &    99& -1&  ~3  &~0 & ~3 &-1 &-1 &-1& -1  &~0 & ~1 & ~1 & -1   &~0   &~0\\
\widehat{\chi}_{13}  &  120  &~0 &-8&  ~3  &~0  &~0  &~0  &~0  &~0 & ~1  &~0  &~0   &~0 & -1  &-1\\
\widehat{\chi}_{14} &  144 & ~4  &~0  &~0 &-3  &~0  &~0& -1 & ~1  &~0  &~0  &~0 & -1&   ~1 &  ~1\\
\widehat{\chi}_{15}&    176 & -4  &~0& -4 &-1  &~0  &~0 & ~1& -1  &~0  &~0  &~0 &  ~1   &~0   &~0\\
\end{smallmatrix}
\right)\label{M12chartable}
\end{equation}
where $\alpha=(\omega +\omega^3 +\omega^4+\omega^5 +\omega^9)$ with $\omega=\exp(2\pi i/11)$.

\subsection{Character Table for $\Mtwo$} \label{M122ch}

The character table for $\Mtwo$ (obtained from the GAP  database\cite{GAP4})
\begin{equation*}
\left(\begin{smallmatrix}
\textrm{Label} & 1A & 2A & 2B & 3A & 3B & 4A & 5A &  6A & 6B & 8A & 10A & 11A & 2C & 4B & 4C & 6C & 10B & 10C & 12A & 12B & 12C\\[3pt]
\widetilde{\chi}_1 & ~1 & ~1 & ~1 & ~1 & ~1 & ~1 & ~1 & ~1 & ~1 & ~1&1&1&~1 & ~1 & ~1 & ~1 & ~1 & ~1 & ~1 & ~1 & ~1 \\
\widetilde{\chi}_2 & ~1 & ~1 & ~1 & ~1 & ~1 & ~1 & ~1 & ~1 & ~1 & ~1 & ~1 & ~1 & -1 & -1 & -1 & -1 & -1 & -1 & -1 & -1 & -1 \\
 \widetilde{\chi}_3 &22 & -2 & ~6 & ~4 & -2 & ~2 & ~2 & -2 & ~0 & ~0 & -2 & ~0 & ~0 & ~0 & ~0 & ~0 & ~0 &
   ~0 & ~0 & ~0 & ~0 \\
\widetilde{\chi}_4 & 32 & ~8 & ~0 & -4 & ~2 & ~0 & ~2 & ~2 & ~0 & ~0 & -2 & -1 & ~0 & ~0 & ~0 & ~0 & ~0 &
   ~0 & ~0 & ~0 & ~0 \\
\widetilde{\chi}_5 & 45 & ~5 & -3 & ~0 & ~3 & ~1 & ~0 & -1 & ~0 & -1 & ~0 & ~1 & ~5 & -3 & ~1 & -1 & ~0
   & ~0 & ~1 & ~0 & ~0 \\
\widetilde{\chi}_6 & 45 & ~5 & -3 & ~0 & ~3 & ~1 & ~0 & -1 & ~0 & -1 & ~0 & ~1 & -5 & ~3 & -1 & ~1 & ~0
   & ~0 & -1 & ~0 & ~0 \\
 \widetilde{\chi}_7 &54 & ~6 & ~6 & ~0 & ~0 & ~2 & -1 & ~0 & ~0 & ~0 & ~1 & -1 & ~0 & ~0 & ~0 & ~0 &
   ~\textrm{A} & -\textrm{A} & ~0 & ~0 & ~0 \\
 \widetilde{\chi}_8 &54 & ~6 &~ 6 & ~0 & ~0 & ~2 & -1 & ~0 & ~0 &~ 0 & ~1 & -1 & ~0 & ~0 & ~0 & ~0 &
   -\textrm{A} & ~\textrm{A} & ~0 & ~0 & ~0 \\
 \widetilde{\chi}_9 &55 & -5 & ~7 & ~1 & 1 & -1 & ~0 & ~1 & ~1 & -1 & ~0 & ~0 & ~5 & ~1 & -1 & -1 & 0
   & ~0 & -1 & ~1 & ~1 \\
 \widetilde{\chi}_{10} &55 & -5 & ~7 & ~1 & 1 & -1 & ~0 & ~1 & ~1 & -1 & ~0 & ~0 & -5 & -1 & ~1 & ~1 & ~0
   & ~0 & ~1 & -1 & -1 \\
 \widetilde{\chi}_{11} & 110 & -10 & -2 & ~2 & ~2 & ~2 & ~0 & ~2 & -2 & ~0 & ~0 & ~0 & ~0 & ~0 & ~0 & ~0 & ~0
   & ~0 & ~0 & ~0 & ~0 \\
 \widetilde{\chi}_{12} &66 & ~6 & ~2 & ~3 & ~0 & -2 & ~1 & ~0 & -1 & ~0 & ~1 & ~0 & ~6 & ~2 & ~0 & ~0 & ~1 & ~1
   & ~0 & -1 & -1 \\
 \widetilde{\chi}_{13}& 66 & ~6 & ~2 & ~3 & ~0 & -2 & ~1 & ~0 & -1 & ~0 & ~1 & ~0 & -6 & -2 & ~0 & ~0 & -1
   & -1 & ~0 & ~1 & ~1 \\
  \widetilde{\chi}_{14}&99 & -1 & ~3 & ~0 & ~3 & -1 & -1 & -1 & ~0 & ~1 & -1 & ~0 & ~1 & -3 & -1 & ~1 &
   ~1 & ~1 & -1 & ~0 & ~0 \\
 \widetilde{\chi}_{15} &99 & -1 & ~3 & ~0 & ~3 & -1 & -1 & -1 & ~0 & ~1 & -1 & ~0 & -1 & ~3 & ~1 & -1 &
   -1 & -1 & ~1 & ~0 &~ 0 \\
 \widetilde{\chi}_{16} &120 & ~0 & -8 & ~3 & ~0 & ~0 & ~0 & ~0 & ~1 & ~0 & ~0 & -1 & ~0 & ~0 & ~0 & ~0 & ~0 &
   ~0 & ~0 & ~\textrm{B} & -\textrm{B} \\
 \widetilde{\chi}_{17} &120 & ~0 & -8 & ~3 & ~0 & ~0 & ~0 & ~0 & ~1 & ~0 & ~0 & -1 & ~0 & ~0 & ~0 & ~0 & ~0 &
   ~0 & ~0 & -\textrm{B} & ~\textrm{B} \\
  \widetilde{\chi}_{18}&144 & ~4 & ~0 & ~0 & -3 & ~0 & -1 & ~1 & ~0 & ~0 & -1 & ~1 & ~4 & ~0 & ~2 & ~1 & -1
   & -1 & -1 & ~0 & ~0 \\
 \widetilde{\chi}_{19} &144 & ~4 & ~0 & ~0 & -3 & ~0 & -1 & ~1 & ~0 & ~0 & -1 & ~1 & -4 & ~0 & -2 & -1 &
   ~1 & ~1 & ~1 & ~0 & ~0 \\
 \widetilde{\chi}_{20} &176 & -4 & ~0 & -4 & -1 & ~0 & ~1 & -1 & ~0 & ~0 & ~1 & ~0 & ~4 & ~0 & -2 & ~1 &
   -1 & -1 & ~1 & ~0 & ~0 \\
 \widetilde{\chi}_{21} &176 & -4 & ~0 & -4 & -1 & ~0 & ~1 & -1 & ~0 & ~0 & ~1 & ~0 & ~4 & ~0 & -2 & ~1 &
   -1 & -1 & ~1 & ~0 & ~0
\end{smallmatrix}
\right)\ ,
\end{equation*}
where $\textrm{A}=(\alpha-\alpha^2-\alpha^3+ \alpha^4)$ and $\textrm{B}=(-\omega^7+\omega^{11})$ with $\alpha=\exp(2\pi i/5)$ and $\omega=\exp(2\pi i/12)$.

\subsection{Decomposing $M_{24}$ characters}\label{decompose}

We decompose some of the $M_{24}$ irreps  that appear in the multiplicative moonshine  in terms of characters of $\Mtwo$ and $M_{12}$. This was obtained using the program GAP\cite{GAP4,Hulpke:2010}:
\begin{align*}
\chi_{23} &= \widetilde{\chi}_2+ \widetilde{\chi}_3 =\widehat{\chi}_{1}+\widehat{\chi}_{2}+\widehat{\chi}_{3}\ ,\\
\chi_{45} &=  \widetilde{\chi}_5 =\widehat{\chi}_{6}\ ,\\
\chi_{231} &=  \widetilde{\chi}_{9}+ \widetilde{\chi}_{20}=\widehat{\chi}_{8}+\widehat{\chi}_{15} \ , \\
\chi_{770} &=  \widetilde{\chi}_{12}+ \widetilde{\chi}_{16}+ \widetilde{\chi}_{17}+ \widetilde{\chi}_{18}+ \widetilde{\chi}_{19}+ \widetilde{\chi}_{20}\\
&=\widehat{\chi}_{11}+2\ \widehat{\chi}_{13}+2\ \widehat{\chi}_{14}+\widehat{\chi}_{15}\ ,
\end{align*}
where the $\Mtwo$ and $M_{12}$ characters are labeled as in the character table given, respectively,  in Appendix \ref{M122ch} and \ref{M12ch}.

\subsection{Multiplicative eta products and $M_{12}$ conjugacy classes}\label{multeta}

Let $\hat{\rho}=1^{a_1}2^{a_2}\cdots N^{a_N}$ be a cycle shape for a partition of $12$. Thus, we have $\sum_i i a_i =12$. Using the map Eq. \eqref{cycleshapemap}, we obtain an eta-product, $g^{\hat{\rho}}(\tau)$. Let $g_{\hat{\rho}}(\tau)$ have the following Fourier expansion
\begin{equation}
g_{\hat{\rho}}(\tau)=\sum_{n=1}^{\infty} a_n\  q^{n/2}\ .
\end{equation}
Extending a definition of Dummit et. al., we call the eta-product multiplicative  if $a_{nm}=a_n a_m$ when gcd$(n,m)=1$. Replacing, $q$ by $q^2$ in the above equation, we see that this reduces precisely to eta-products considered by Dummit et. al.\cite{Dummit:1985}.  They found 30 multiplicative eta-products -- among these eta-products, we need to look at cycles shapes with only \textit{even} cycles so that it can be reduced to a cycle shape with $\sum_i i a_i =12$. We find 15 of the 30 cycle shapes in their list satisfy our criterion. Further, 12 of the 15 cycle shapes also arise as $M_{12}$ conjugacy classes. The cycles shapes $4^3$, $3\ 9$ and $12$ gives rise to multiplicative eta-products but are \textit{not} $M_{12}$ conjugacy classes. The $M_{12}$ classes $4B$ and $8B$ do not give rise to multiplicative eta-products. We list them in the Table  \ref{M12balanced}.
\begin{table}[hbt]
\centering
\newcommand\T{\rule{0pt}{2.6ex}}
\begin{tabular}{cccccccccccc} \hline
 $1^{12}$ & $2^6$ & $1^42^4$ & $1^3 3^3$ & $3^4$ & $2^24^2$ & $1^2 5^2$ & $6^2$ & $1\ 2\ 3\ 6$ & $4\ 8$ &  $2\ 10$ & $1\ 11$\T \\[2pt] \hline
 1A & 2A & 2B & 3A & 3B & 4A & 5A & 6A &6B & 8A &  10A & 11A/B \T \\[2pt] \hline
\end{tabular}
\caption{Balanced cycles shapes associated with multiplicative eta-products and their $M_{12}$ conjugacy class. }\label{M12balanced}
\end{table}

\section{The modular forms}

The genus-one theta functions are defined by
\begin{equation}
\theta\left[\genfrac{}{}{0pt}{}{a}{b}\right] \left(z_1,z_2\right)
=\sum_{l \in \BZ} 
q^{\frac12 (l+\frac{a}2)^2}\ 
r^{(l+\frac{a}2)}\ e^{i\pi lb}\ ,
\end{equation}
where $a.b\in (0,1)\mod 2$ and $q=\exp(2\pi i z_1)$ and $r=\exp(2\pi i z_2)$.
 One has $\vartheta_1 
\left(z_1,z_2\right)\equiv\theta\left[\genfrac{}{}{0pt}{}{1}{1}\right] 
\left(z_1,z_2\right)$, $\vartheta_2 
\left(z_1,z_2\right)\equiv\theta\left[\genfrac{}{}{0pt}{}{1}{0}\right] 
\left(z_1,z_2\right)$, $\vartheta_3 
\left(z_1,z_2\right)\equiv\theta\left[\genfrac{}{}{0pt}{}{0}{0}\right] 
\left(z_1,z_2\right)$ and $\vartheta_4 
\left(z_1,z_2\right)\equiv\theta\left[\genfrac{}{}{0pt}{}{0}{1}\right] 
\left(z_1,z_2\right)$.

We define the genus-two theta constants as follows\cite{Nikulin:1995}:
\begin{equation}
\theta\left[\genfrac{}{}{0pt}{}{\mathbf{a}}{\mathbf{b}}\right]
\left(\mathbf{Z}\right)
=\sum_{(l_1, l_2)\in \BZ^2} 
q^{\frac12 (l_1+\frac{a_1}2)^2}\ 
r^{(l_1+\frac{a_1}2)(l_2+\frac{a_2}2)}\ 
s^{\frac12 (l_2+\frac{a_2}2)^2}\ 
e^{i\pi(l_1b_1+l_2b_2)}\ ,
\end{equation}
where $\mathbf{a}=\begin{pmatrix}a_1\\ a_2
\end{pmatrix}$,
$\mathbf{b}=\begin{pmatrix}b_1\\ b_2
\end{pmatrix}$,
and $\mathbf{Z}=\begin{pmatrix}z_1 & z_2 \\ z_2 &
z_3\end{pmatrix}\in \mathbb{H}_2$. 
Further, we have defined $q=\exp(2\pi i z_1)$,
$r=\exp(2\pi i z_2)$ and $s=\exp(2\pi i z_3)$.
The constants $(a_1,a_2,b_1,b_2)$ take values $(0,1)$. 
Thus there are sixteen genus-two theta constants. The
even theta constants are those for which
$\mathbf{a}^{\textrm{T}}\mathbf{b}=0\mod 2$. There are ten such theta
constants.  Note that six of the even theta constants with $\mathbf{a}\neq0$ have even Fourier coefficients while the remaining four theta constants with $\mathbf{a}=0$ have integral Fourier coefficients.

\begin{equation}
\Phi^{4B}_3(\mathbf{Z})=\left(\frac18\ \theta\!\left[\begin{smallmatrix} 1\\ 0 \\ 0 \\ 1\end{smallmatrix}\right]\!\!\left(2\mathbf{Z}\right)\ 
\theta\!\left[\begin{smallmatrix} 0\\ 1 \\ 1 \\ 0\end{smallmatrix}\right]\!\!\left(2\mathbf{Z}\right)\ 
\theta\!\left[\begin{smallmatrix} 1\\ 1 \\ 1 \\ 1\end{smallmatrix}\right]\!\!\left(2\mathbf{Z}\right)\ 
\right)^2\equiv\left[\Delta_{3/2}(\mathbf{Z})\right]^2\ .
\end{equation}

\begin{equation}
\Delta^{4B}_1(\mathbf{Z})  = \frac14\ \theta\!\left[\begin{smallmatrix} 0\\ 1 \\ 1 \\ 0\end{smallmatrix}\right]\!\!\left(2\mathbf{Z}\right)\ 
\theta\!\left[\begin{smallmatrix} 1\\ 1 \\ 1 \\ 1\end{smallmatrix}\right]\!\!\left(2\mathbf{Z}\right)\ .
\end{equation}
\bibliography{master}
\end{document}